\documentclass[11pt]{article}
\usepackage{amssymb,latexsym}
\usepackage{epsfig}
\usepackage{eufrak}
\usepackage{amsmath}
\usepackage{mathrsfs}
\usepackage{color}

\newtheorem{theorem}{Theorem}
\newtheorem{lemma}{Lemma}

\newcommand{\be}{\begin{equation}}
\newcommand{\ee}{\end{equation}}
\newcommand{\bee}{\begin{eqnarray*}}
\newcommand{\eee}{\end{eqnarray*}}
\newcommand{\bel}{\begin{eqnarray}}
\newcommand{\eel}{\end{eqnarray}}
\newcommand{\bec}{\begin{cases}}
\newcommand{\eec}{\end{cases}}
\newcommand{\bem}{\begin{bmatrix}}
\newcommand{\eem}{\end{bmatrix}}

\newcommand{\la}{\label}
\newcommand{\li}{\left}
\newcommand{\ri}{\right}

\newcommand{\lc}{\lceil}
\newcommand{\rc}{\rceil}
\newcommand{\lf}{\lfloor}
\newcommand{\rf}{\rfloor}

\newcommand{\vep}{\varepsilon}
\newcommand{\lm}{\lambda}

\newcommand{\de}{\delta}

\newcommand{\se}{\theta}

\newcommand{\al}{\alpha}
\newcommand{\ba}{\beta}

\newcommand{\Om}{\Omega}

\newcommand{\f}{\frac}

\newcommand{\cd}{\cdots}

\newcommand{\qqu}{\qquad}

\newcommand{\mscr}{\mathscr}

\newcommand{\bb}{\mathbb}

\newcommand{\wh}{\widehat}

\newcommand{\mrm}{\mathrm}
\newcommand{\bs}{\boldsymbol}

\newcommand{\tx}{\text}

\newcommand{\pa}{\partial}

\newcommand{\bed}{\begin{description}}
\newcommand{\eed}{\end{description}}
\newcommand{\bei}{\begin{itemize}}
\newcommand{\eei}{\end{itemize}}
\newcommand{\ben}{\begin{enumerate}}
\newcommand{\een}{\end{enumerate}}

\newcommand{\beL}{\begin{lemma}}
\newcommand{\eeL}{\end{lemma}}
\newcommand{\beT}{\begin{theorem}}
\newcommand{\eeT}{\end{theorem}}
\newcommand{\sect}{\section}

\newcommand{\bpf}{\begin{pf}}
\newcommand{\epf}{\end{pf}}
\newcommand{\bsk}{\bigskip}

\newcommand{\pfbox}{\hfill\mbox{$\Box$}}

\newenvironment{pf}{\paragraph*{Proof{\rm.}}}{\pfbox\bigskip}

\begin{document}

\title{{\bf Exact Computation of Minimum Sample size for Estimation of Poisson Parameters}
\thanks{The author had been previously working with Louisiana State University at Baton Rouge, LA 70803, USA,  and is now with
Department of Electrical Engineering, Southern University and A\&M College, Baton Rouge, LA 70813, USA; Email:
chenxinjia@gmail.com}}

\author{Xinjia Chen}

\date{July 2007}

\maketitle

\begin{abstract}

In this paper, we develop an approach for the exact determination of the minimum sample size for the estimation
of a Poisson parameter with prescribed margin of error and confidence level.  The exact computation is made
possible by reducing infinite many evaluations of coverage probability to finite many evaluations.  Such
reduction is based on our discovery that the minimum of coverage probability with respect to a Poisson parameter
bounded in an interval is attained at a discrete set of finite many values.

\end{abstract}

\sect{Introduction} The estimation of a Poisson parameter finds numerous applications in various fields of
sciences and engineering \cite{Desu}.  The problem is formulated as follows.

Let $X$ be a Poisson random variable defined in a probability space $(\Om, \mscr{F}, \Pr)$ such that $\Pr \{ X =
k \} =  \f{ \lm^k e^{- \lm} } { k!  }, \; k = 0, 1, \cd$, where $\lm > 0$ is referred to as a Poisson parameter.
It is a frequent problem to estimate $\lm$ based on $n$ identical and independent samples $X_1, \cd, X_n$ of
$X$.

An estimate of $\lm$ is conventionally taken as $\wh{\bs{\lm}}_n = \f{ \sum_{i=1}^n X_i } { n }$.  The nice
property of such estimate is that it is of maximum likely-hood and possesses minimum variance among all unbiased
estimates.  A crucial question in the estimation is as follows:

\bsk

{\it Given the knowledge that $\lm$ belongs to interval $[a, b]$, what is the minimum sample size $n$ that
guarantees the difference between $\wh{\bs{\lm}}_n$ and $\lm$ be bounded within some prescribed margin of error
with a confidence level higher than a prescribed value?}

\bsk

The main contribution of this paper is to provide exact answer to this important question. The paper is
organized as follows. In Section 2, the techniques for computing the minimum sample size is developed with the
margin of error taken as a bound of absolute error. In Section 3, we derive corresponding sample size method by
using relative error bound as the margin of error.  In Section 4, we develop techniques for computing minimum
sample size with a mixed error criterion.  Section 5 is the conclusion. The proofs are given in Appendices.

Throughout this paper, we shall use the following notations. The set of integers is denoted by $\bb{Z}$.  The
ceiling function and floor function are denoted respectively by $\lc . \rc$ and $\lf . \rf$ (i.e., $\lc x \rc$
represents the smallest integer no less than $x$; $\lf x \rf$ represents the largest integer no greater than
$x$). The multivariate function $S(n, k, l, \lm)$ means {\small $S(n, k, l, \lm) =  \sum_{i = k}^l \f{ \lm^i
e^{- \lm} }{ i! }$}. The left limit as $\eta$ tends to $0$ is denoted as $\lim_{\eta \downarrow 0}$.  The other
notations will be made clear as we proceed.

\section{Control of Absolute Error}

Let $\vep \in (0,1)$ be the margin of absolute error and $\de \in
(0,1)$ be the confidence parameter.  In many applications, it is
desirable to find the minimum sample size $n$ such that \[
 \Pr \li
\{ |\wh{\bs{\lm}}_n - \lm| < \vep \ri \} > 1 - \de \] for any $\lm \in [a,b]$. Here $\Pr \li \{ |\wh{\bs{\lm}}_n
- \lm | < \vep \ri \}$ is referred to as the coverage probability.  The interval $[a, b]$ is introduced to take
into account the knowledge of $\lm$.  The exact determination of minimum sample size is readily tractable with
modern computational power by taking advantage of the behavior of the coverage probability characterized by
Theorem \ref{thm_abs} as follows.

\bsk

\beT \la{thm_abs} Let $0 < \vep < 1$ and $0 \leq a < b$. Let $X_1,
\cd, X_n$ be identical and independent Poisson random variables
with mean $\lm \in [a, b]$. Let $ \wh{\bs{\lm}}_n = \f{ \sum_{i =
1}^n X_i  }{ n }$. Then, the minimum of $\Pr \{ | \wh{\bs{\lm}}_n
- \lm | < \vep \}$ with respect to $\lm \in [a, b]$ is achieved at
the finite set $\{a, b \} \cup
 \{ \f{\ell}{n} + \vep \in (a, b) : \ell \in \bb{Z} \} \cup
 \{ \f{\ell}{n} - \vep \in (a, b) : \ell \in \bb{Z} \} $,
 which has less than $2n(b-a) + 4$ elements.
\eeT

See Appendix A for a proof.  The application of Theorem \ref{thm_abs} in the computation of minimum sample size
is obvious.  For a fixed sample size $n$, since the minimum of coverage probability with $\lm \in [a, b]$ is
attained at a finite set, it can determined by a computer whether the sample size $n$ is large enough to ensure
$\Pr \li \{ |\wh{\bs{\lm}}_n - \lm| < \vep \ri \} > 1 - \de$ for any $\lm \in [a, b]$. Starting from $n = 2$,
one can find the minimum sample size by gradually incrementing $n$ and checking whether $n$ is large enough.

\section{Control of Relative Error}

Let $\vep \in (0,1)$ be the margin of  relative error and $\de \in
(0,1)$ be the confidence parameter.  It is interesting to determine
the minimum sample size $n$ so that
\[
\Pr \li \{  \li | \f{ \wh{\bs{\lm}}_n - \lm } {\lm} \ri | < \vep
\ri \}
> 1 - \de
\]
for any $\lm \in [a,b]$.  As has been pointed out in Section 2, an essential machinery is to reduce infinite
many evaluations of the coverage probability $\Pr \{ | \wh{\bs{\lm}}_n - \lm | < \vep \lm \}$ to finite many
evaluations. Such reduction can be accomplished by making use of Theorem \ref{thm_rev} as follows.

 \beT \la{thm_rev} Let $0 < \vep < 1$ and $0
< a  < b$.  Let $X_1, \cd, X_n$ be identical and independent
Poisson random variables  with mean $\lm \in [a, b]$. Let $
\wh{\bs{\lm}}_n = \f{ \sum_{i = 1}^n X_i  }{ n }$. Then, the
minimum of $\Pr \li \{  \f{ |\wh{\bs{\lm}}_n - \lm| } {\lm} < \vep
\ri \}$ with respect to $\lm \in [a, b]$ is achieved at the finite
set $\{a, b \} \cup
 \{ \f{\ell}{n(1 + \vep)} \in (a, b) : \ell \in \bb{Z} \} \cup
 \{ \f{\ell}{n (1 - \vep)} \in (a, b) : \ell \in \bb{Z} \}$,
 which has less than $2n(b-a) + 4$ elements.
\eeT

See Appendix B for a proof.

\section{Control of Absolute Error or Relative Error}

Let $\vep_a \in (0,1)$ and $\vep_r \in (0,1)$ be respectively the
margins of absolute error and relative error.  Let $\de \in (0,1)$
be the confidence parameter.  In many situations, it is desirable to
find the smallest sample size $n$ such that \be \la{rule2} \Pr \li
\{ |\wh{\bs{\lm}}_n - \lm| < \vep_a \; \; \mrm{or} \;\;  \li | \f{
\wh{\bs{\lm}}_n - \lm } {\lm} \ri | < \vep_r \ri \}
> 1 - \de \ee for any $\lm \in
[a,b]$.  To make it possible to compute exactly
 the minimum sample size associated with (\ref{rule2}), we have
Theorem \ref{thm_abs_rev} as follows.

\beT \la{thm_abs_rev} Let $0 < \vep_a < 1, \; 0 < \vep_r < 1$ and
$0 \leq a < \f{\vep_a}{\vep_r} < b$.
 Let $X_1, \cd, X_n$ be identical and independent
Poisson random variables with mean $\lm \in [a, b]$. Let $
\wh{\bs{\lm}}_n = \f{ \sum_{i = 1}^n X_i  }{ n }$. Then, the minimum
of $\Pr \li \{ |\wh{\bs{\lm}}_n - \lm| < \vep_a \; \; \mrm{or} \;\;
\li | \f{ \wh{\bs{\lm}}_n - \lm } {\lm} \ri | < \vep_r  \ri \}$ with
respect to $\lm \in [a, b]$ is achieved at the finite set $\{a, b,
\f{\vep_a}{\vep_r} \} \cup
 \{ \f{\ell}{n} + \vep_a \in (a, \f{\vep_a}{\vep_r}) : \ell \in \bb{Z} \} \cup
 \{ \f{\ell}{n} - \vep_a \in (\f{\vep_a}{\vep_r}, b) : \ell \in \bb{Z} \} \cup
 \{ \f{\ell}{n(1 + \vep_r)} \in (a, \f{\vep_a}{\vep_r}) : \ell \in \bb{Z} \} \cup
 \{ \f{\ell}{n (1 - \vep_r)} \in (\f{\vep_a}{\vep_r}, b) : \ell \in \bb{Z}
 \}$, which has less than $2n(b-a) + 7$ elements.
\eeT

\bsk

Theorem \ref{thm_abs_rev} can be shown by applying Theorem
\ref{thm_abs} and Theorem \ref{thm_rev} with the observation that
\[
\Pr \li \{ |\wh{\bs{\lm}}_n - \lm| < \vep_a \; \; \mrm{or} \;\; \li
| \f{ \wh{\bs{\lm}}_n - \lm } {\lm} \ri | < \vep_r \ri \} = \bec \Pr
\li \{ |\wh{\bs{\lm}}_n - \lm| < \vep_a  \ri \} & \tx{for} \; \lm
\in \li [ a,
\f{\vep_a}{\vep_r} \ri ], \\
\Pr \li \{ \li | \f{ \wh{\bs{\lm}}_n - \lm } {\lm} \ri | < \vep_r
\ri \} & \tx{for} \; \lm \in \li ( \f{\vep_a}{\vep_r}, b \ri ]. \eec
\]

\bsk

By virtue of Chernoff bounds, it can be shown that, for any $\vep
\in (0,1)$,
\[
\Pr \{ \wh{\bs{\lm}}_n \leq (1 - \vep) \lm \} < \li [ \f{e^{- \vep}}
{ (1-\vep)^{1-\vep}  } \ri ]^{n \lm} < \exp \li ( - \f{\lm n
\vep^2}{2} \ri ), \] \[
 \Pr \{ \wh{\bs{\lm}}_n \geq (1 + \vep) \lm
\} < \li [ \f{e^{ \vep}} { (1 + \vep)^{1 + \vep}  } \ri ]^{n \lm} <
\exp \li ( - (2 \ln 2 - 1) \lm n \vep^2 \ri ).
\]

As a result, $\Pr \{ | \wh{\bs{\lm}}_n - \lm | > \vep \lm \} < \de$ if
\[
\lm > \f{ \ln \f{2}{\de}  } { (2 \ln 2 -1) n \vep^2  }.
\]
Therefore, to check whether (\ref{rule2}) is satisfied for any $\lm \in [a, b]$, it suffices to check
(\ref{rule2}) for
\[
a \leq \lm \leq \min \li \{ b, \; \f{ \ln \f{2}{\de}  } { (2 \ln 2 -1) n \vep_r^2  } \ri \}.
\]

Finally, we would like to point out that similar characteristics of the coverage probability can be shown for
the problem of estimating binomial parameter or the proportion of finite population, which allows for the exact
computation of minimum sample size. For details, see our recent papers \cite{Chen, Chen2}.

\section{Conclusion}

We have developed an exact method for the computation of minimum sample size for the estimation of Poisson
parameters, which only requires finite many evaluations of the coverage probability. Our sample size method
permits rigorous control of statistical sampling error.

\appendix

\sect{Proof of Theorem \ref{thm_abs}}

Define $K = \sum_{i=1}^n X_i$ and \[ C(\lm)  =  \Pr \li \{ \li |
\f{K}{n} - \lm \ri | < \vep \ri \} = \Pr \li \{ g(\lm) \leq K \leq
h(\lm) \ri \}
\] where
\[
g(\lm) = \max ( 0, \lf n( \lm - \vep) \rf + 1 ), \qqu h(\lm) = \lc
n( \lm + \vep) \rc - 1.
\]
It should be noted that $C(\lm), \; g(\lm)$ and $h(\lm)$ are actually multivariate functions of $\lm, \; \vep$
and $n$.  For simplicity of notations, we drop the arguments $n$ and $\vep$ throughout the proof of Theorem
\ref{thm_abs}.

We need some preliminary results.

\beL \la{minus} Let $\lm_\ell = \f{\ell}{n} - \vep$ where $\ell
\in \bb{Z}$. Then, $h(\lm) = h(\lm_{\ell + 1}) = \ell$ for any
$\lm \in (\lm_\ell, \lm_{\ell +1})$. \eeL

\bpf For $\lm \in ( \lm_\ell, \; \lm_{\ell + 1})$, we have $0 < n
\li (\lm - \lm_\ell \ri ) < 1$ and \bee h (\lm)  & = &
\lc n( \lm + \vep) \rc - 1\\
& = &  \lc n( \lm_\ell + \vep + \lm - \lm_\ell ) \rc - 1\\
& = & \li \lc n \li ( \f{\ell}{n} - \vep  + \vep + \lm - \lm_\ell  \ri )  \ri \rc - 1\\
& = & \ell - 1 + \li \lc n \li (\lm - \lm_\ell \ri )  \ri \rc\\
& = & \ell\\
& = & \li \lc n \li ( \f{\ell + 1}{n} - \vep  + \vep \ri )  \ri \rc
- 1 = h(\lm_{\ell + 1}). \eee

\epf

\beL \la{plus} Let $\lm_\ell = \f{\ell}{n} + \vep$ where $\ell \in
\bb{Z}$. Then, $g(\lm) = g(\lm_{\ell})= \max \{0, \ell + 1 \}$ for
any $\lm \in (\lm_\ell, \lm_{\ell +1})$. \eeL

 \bpf
For $\lm \in \li ( \lm_\ell, \; \lm_{\ell + 1} \ri )$,  we have
$-1 < n \li (\lm - \lm_{\ell + 1} \ri ) < 0$ and \bee
g (\lm) & = & \max ( 0, \lf n( \lm - \vep) \rf + 1 )\\
& = &  \max ( 0, \lf n( \lm_{\ell + 1} - \vep + \lm - \lm_{\ell + 1} ) \rf + 1 )\\
& = &  \max \li ( 0, \li \lf n \li ( \f{ \ell + 1 } { n  } + \vep -
\vep \ri ) \ri \rf +
\lf n( \lm - \lm_{\ell + 1} ) \rf +  1 \ri )\\
& = & \max \li ( 0, \li \lf n \li ( \f{ \ell + 1 } { n  } + \vep  - \vep \ri ) \ri \rf - 1 +  1 \ri )\\
& = & \max \{0, \ell + 1 \}\\
& = & \max \li ( 0, \li \lf n \li ( \f{ \ell } { n  } + \vep - \vep
\ri ) \ri \rf +  1 \ri ) =  g(\lm_\ell). \eee

\epf

\beL \la{constant} Let $\al < \ba$ be two consecutive elements of
the ascending arrangement of all distinct elements of $\{a, b \}
\cup
 \{ \f{\ell}{n} + \vep \in (a, b) : \ell \in \bb{Z} \} \cup
 \{ \f{\ell}{n} - \vep \in (a, b) : \ell \in \bb{Z} \} $.
Then, both $g(\lm)$ and $h (\lm)$ are constants for any $\lm \in (\al, \ba)$.
 \eeL

 \bpf
Since $\al$ and $\ba$ are two consecutive elements of the
ascending arrangement of all distinct elements of the set, it must
be true that there is no integer $\ell$ such that {\small $\al <
\f{\ell}{n} + \vep < \ba$} or {\small $\al < \f{\ell}{n} - \vep <
\ba$}.  It follows that there exist two integers $\ell$ and
$\ell^\prime$ such that {\small $(\al, \ba) \subseteq \li (
\f{\ell}{n} + \vep, \f{\ell + 1}{n} + \vep \ri )$} and {\small
$(\al, \ba) \subseteq \li ( \f{\ell^\prime}{n} - \vep,
\f{\ell^\prime + 1}{n} - \vep \ri )$}. Applying Lemma \ref{minus}
and Lemma \ref{plus}, we have {\small $g(\lm) = g \li (
\f{\ell}{n} + \vep \ri )$} and {\small $h(\lm) = h \li
(\f{\ell^\prime + 1}{n} - \vep \ri )$} for any $\lm \in (\al,
\ba)$.

 \epf

\beL \la{lem_lim}
 For any $\lm \in (0,1)$, $\lim_{\eta \downarrow 0} C(\lm + \eta) \geq C(\lm)$
 and $\lim_{\eta \downarrow 0} C(\lm - \eta) \geq C(\lm)$.
\eeL

\bpf

Observing that $h(\lm + \eta) \geq h(\lm)$ for any $\eta > 0$ and
that
 \bee g(\lm + \eta) & =  & \max ( 0, \lf n( \lm + \eta - \vep)
\rf + 1 ) \\
& = & \max ( 0, \lf n( \lm - \vep) \rf + 1 + \lf n( \lm - \vep) -
\lf n( \lm  - \vep) \rf + n \eta \rf )\\
 & = & \max ( 0, \lf n( \lm - \vep)
\rf + 1 )  = g(\lm) \eee for $0 < \eta < \f{ 1 + \lf n( \lm - \vep)
\rf -  n( \lm - \vep)} {n}$, we have \be \la{ineqa} S(n, g(\lm +
\eta), h (\lm + \eta), \lm + \eta ) \geq S(n, g(\lm), h (\lm), \lm +
\eta ) \ee for $0 < \eta < \f{ 1 + \lf n( \lm - \vep) \rf - n( \lm -
\vep)} {n}$.  Since
 \[ h(\lm + \eta)  =  \lc n( \lm + \eta + \vep)
\rc - 1 = \lc n( \lm + \vep) \rc - 1 + \lc n( \lm + \vep) - \lc n(
\lm + \vep) \rc + n \eta \rc, \] we have
\[
h(\lm + \eta) = \bec \lc n( \lm + \vep) \rc & \tx{for} \; n( \lm +
\vep) = \lc n( \lm  + \vep) \rc \; \tx{and} \; 0 < \eta <
\f{1}{n},\\
\lc n( \lm + \vep) \rc  - 1 & \tx{for} \; n( \lm + \vep) \neq \lc n(
\lm + \vep) \rc\; \tx{and} \; 0 < \eta < \f{\lc n( \lm + \vep) \rc -
n( \lm + \vep)}{n}. \eec
\]
It follows that both $g(\lm + \eta)$ and $h(\lm + \eta)$ are
independent of $\eta$ if $\eta > 0$ is small enough.  Since $S(n, g,
h, \lm + \eta)$ is continuous with respect to $\eta$ for fixed $g$
and $h$, we have that $\lim_{\eta \downarrow 0} S(n, g(\lm + \eta),
h (\lm + \eta), \lm + \eta )$ exists.  As a result, \bee \lim_{\eta
\downarrow 0} C(\lm + \eta) & = & \lim_{\eta \downarrow 0} S(n,
g(\lm + \eta), h (\lm + \eta), \lm + \eta )\\
& \geq & \lim_{\eta \downarrow 0} S(n, g(\lm), h (\lm), \lm + \eta )
= S(n,  g(\lm), h (\lm), \lm ) = C(\lm), \eee where the inequality
follows from (\ref{ineqa}).

Observing that $g(\lm - \eta) \leq g(\lm)$ for any $\eta > 0$ and
that
 \bee h(\lm - \eta) & =  & \lc n( \lm - \eta + \vep)
\rc - 1\\
& = & \lc n( \lm + \vep) \rc - 1 + \lc n( \lm + \vep) -
\lc n( \lm  + \vep) \rc - n \eta \rc\\
 & = & \lc n( \lm + \vep)
\rc - 1   = h(\lm) \eee for $0 < \eta < \f{ 1 + n( \lm + \vep) - \lc
n( \lm + \vep) \rc } {n}$, we have \be \la{ineqb} S(n, g(\lm -
\eta), h (\lm - \eta), \lm - \eta ) \geq S(n, g(\lm), h (\lm), \lm -
\eta ) \ee for {\small $0 < \eta < \min \li \{ \lm, \f{ 1 + n( \lm +
\vep) - \lc n( \lm + \vep) \rc } {n} \ri \}$}.  Since \bee g(\lm -
\eta) & = & \max ( 0, \lf n( \lm - \eta - \vep)
\rf + 1 ) \\
& = & \max ( 0, \lf n( \lm - \vep) \rf + 1 + \lf n( \lm - \vep) -
\lf n( \lm  - \vep) \rf - n \eta \rf ), \eee we have {\small \[
g(\lm - \eta) = \bec \max ( 0, \lf n( \lm - \vep) \rf ) & \tx{for}
\; n( \lm - \vep) = \lf n( \lm  - \vep) \rf \; \tx{and} \; 0 < \eta
< \f{1}{n},\\
\max ( 0, \lf n( \lm - \vep) \rf  + 1) & \tx{for} \; n( \lm -
\vep) \neq \lf n( \lm - \vep) \rf \; \tx{and} \; 0 < \eta < \f{n(
\lm - \vep) - \lf n( \lm - \vep) \rf }{n}. \eec
\]}
It follows that both $g(\lm - \eta)$ and $h(\lm - \eta)$ are
independent of $\eta$ if $\eta > 0$ is small enough.  Since $S(n, g,
h, \lm - \eta)$ is continuous with respect to $\eta$ for fixed $g$
and $h$, we have that $\lim_{\eta \downarrow 0} S(n, g(\lm - \eta),
h (\lm - \eta), \lm - \eta )$ exists. Hence, \bee \lim_{\eta
\downarrow 0} C(\lm - \eta) & =
& \lim_{\eta \downarrow 0} S(n, g(\lm - \eta), h (\lm - \eta), \lm -\eta )\\
& \geq & \lim_{\eta \downarrow 0} S(n, g(\lm), h (\lm), \lm -\eta )
=  S(n, g(\lm), h (\lm), \lm ) = C(\lm), \eee where the inequality
follows from (\ref{ineqb}).

\epf

\beL \la{inbetween} Let $\al < \ba$ be two consecutive elements of
the ascending arrangement of all distinct elements of $\{a, b \}
\cup
 \{ \f{\ell}{n} + \vep \in (a, b) : \ell \in \bb{Z} \} \cup
 \{ \f{\ell}{n} - \vep \in (a, b) : \ell \in \bb{Z} \} $.  Then,
 $C(\lm) \geq \min \{ C(\al), \; C(\ba) \}$ for any $\lm \in (\al, \ba)$.
\eeL

\bpf

By Lemma \ref{constant}, both $g(\lm)$ and $h(\lm)$ are constants for any $\lm \in (\al, \ba)$. Hence, we can
drop the argument and write $g(\lm) = g, \; h(\lm) = h$ and $C(\lm) = S(n, g, h, \lm)$.

For $\lm \in (\al, \ba)$, define interval $[\al + \eta, \ba - \eta]$
with {\small $0 < \eta < \min \li ( \lm - \al, \ba - \lm, \f{\ba -
\al} { 2 } \ri )$}. Then, $C(\lm) \geq \min_{\mu \in [\al + \eta,
\ba - \eta] } C(\mu)$. Note that $\f{ \pa S(n, 0, l, \lm) } { \pa
\lm } = - \f{\lm^l e^{-\lm} }{l!}$ and thus, for $g > 0$, {\small
\bee \f{ \pa S(n, g, h, \lm) } { \pa \lm }
& = & \f{ \pa S(n, 0, h, \lm) } { \pa \lm } - \f{ \pa S(n, 0, g-1, \lm) } { \pa \lm }\\
& = & \f{\lm^{g-1} e^{-\lm} }{(g-1)!} - \f{\lm^h e^{-\lm} }{h!}\\
& = & \li [ \f{ h !   } { (g-1) ! } - \lm^ {h-g+1} \ri ] \f{
\lm^{g-1} e^{-\lm} } { h !}  > 0 \eee} if {\small $\lm <  \li [ \f{
h ! } { (g-1) !
 } \ri ]^{\f{1}{h-g+1}}$. } From such investigation of the derivative of $S(n, g, h, \lm)$ with
respective to $\lm$, we can see that, for {\small $0 < \eta < \min \li ( \lm - \al, \ba - \lm, \f{\ba - \al} { 2
} \ri )$}, one of the following three cases must be true: (1) $C(\mu)$ decreases monotonically for $\mu \in [\al
+ \eta, \ba - \eta]$; (2) $C(\mu)$ increases monotonically for $\mu \in [\al + \eta, \ba - \eta]$; (3) there
exists a number $\se \in (\al + \eta, \ba - \eta)$ such that $C(\mu)$ increases monotonically for $\mu \in [\al
+ \eta, \se]$ and decreases monotonically for $\mu \in (\se, \ba - \eta]$. It follows that
\[
C(\lm) \geq \min_{\mu \in [\al + \eta, \ba - \eta]} C(\mu) = \min \{
C(\al + \eta), \; C(\ba - \eta) \}
\]
for {\small $0 < \eta < \min \li ( \lm - \al, \ba - \lm, \f{\ba -
\al} { 2 } \ri )$}. By Lemma \ref{lem_lim}, both $\lim_{\eta
\downarrow 0} C(\al + \eta)$ and $\lim_{\eta \downarrow 0}C(\ba -
\eta)$ exist and \bee C(\lm) & \geq & \lim_{\eta \downarrow 0} \;
\min \{
C(\al + \eta), \; C(\ba - \eta) \}\\
& = & \min \li \{ \lim_{\eta \downarrow 0} C(\al + \eta), \;
\lim_{\eta \downarrow 0} C(\ba - \eta) \ri \} \geq \min \{ C(\al),
\; C(\ba) \} \eee for any $\lm \in (\al, \ba)$. \epf

\bsk

Finally, to show Theorem \ref{thm_abs}, note that the statement about the coverage probability  follows
immediately from Lemma \ref{inbetween}.  The number of elements of the finite set can be calculated by using the
property of the ceiling and floor functions.

\sect{Proof of Theorem \ref{thm_rev} }

Define \[ C(\lm)  =  \Pr \li \{ \li | \f{K}{n} - \lm \ri | < \vep
\lm \ri \} = \Pr \li \{ g(\lm) \leq K \leq h(\lm) \ri \}
\] where
\[
g(\lm) = \lf n \lm (1 - \vep) \rf + 1 , \qqu h(\lm) = \lc n \lm (1 +
\vep ) \rc - 1.
\]
It should be noted that $C(\lm), \; g(\lm)$ and $h(\lm)$ are actually multivariate functions of $\lm, \; \vep$
and $n$.  For simplicity of notations, we drop the arguments $n$ and $\vep$ throughout the proof of Theorem
\ref{thm_rev}.

We need some preliminary results.

\beL \la{minus_rev} Let $\lm_\ell = \f{\ell}{n (1 + \vep)}$ where
$\ell \in \bb{Z}$. Then, $h(\lm) = h(\lm_{\ell + 1}) = \ell$ for
any $\lm \in (\lm_\ell, \lm_{\ell +1})$. \eeL

\bpf For $\lm \in ( \lm_\ell, \; \lm_{\ell + 1})$, we have $0 < n
(1 + \vep) \li (\lm - \lm_\ell \ri ) < 1$ and \bee h (\lm)  & = &
\lc n \lm (1 + \vep) \rc - 1\\
& = & \lc n \lm_\ell (1 + \vep) + (1 + \vep) ( \lm - \lm_\ell ) \rc - 1\\
& = & \li \lc n \li [ \f{\ell}{n}  + (1 + \vep) ( \lm - \lm_\ell)  \ri ]  \ri \rc - 1\\
& = & \ell - 1 + \li \lc n (1 +
\vep) \li (\lm - \lm_\ell \ri )  \ri \rc\\
& = & \ell\\
& = & \li \lc n \li [ \f{\ell + 1}{n (1 + \vep)}  \times (1 + \vep)
\ri ]  \ri \rc - 1 = h(\lm_{\ell + 1}). \eee

\epf

\beL \la{plus_rev} Let $\lm_\ell = \f{\ell}{n (1 - \vep)}$ where
$\ell \in \bb{Z}$. Then, $g(\lm) = g(\lm_{\ell})= \ell + 1$ for
any $\lm \in (\lm_\ell, \lm_{\ell +1})$. \eeL

 \bpf
For $\lm \in \li ( \lm_\ell, \; \lm_{\ell + 1} \ri )$,  we have
$-1 < n (1 - \vep)  \li (\lm - \lm_{\ell + 1} \ri ) < 0$ and \bee
g (\lm) & = &  \lf n \lm (1 - \vep) \rf + 1 \\
& = & \lf n[ \lm_{\ell + 1} (1 - \vep) + (1 - \vep) (\lm - \lm_{\ell + 1}) ] \rf + 1\\
& = &  \li \lf n \times \f{ \ell + 1 } { n (1 - \vep) } \times (1 -
\vep)  \ri \rf +
\lf n (1 - \vep) ( \lm - \lm_{\ell + 1} ) \rf +  1 \\
& = & \li \lf n \times \f{ \ell + 1 } { n (1 - \vep) } \times (1 -
\vep)  \ri \rf   - 1 +  1 \\
& = & \ell + 1\\
& = & \li \lf n \times \f{ \ell } { n (1 - \vep) } \times (1 - \vep)
\ri \rf +  1 = g(\lm_\ell). \eee

\epf

\beL \la{constant_rev} Let $\al < \ba$ be two consecutive elements
of the ascending arrangement of all distinct elements of $\{a, b
\} \cup
 \{ \f{\ell}{n (1 - \vep)}  \in (a, b) : \ell \in \bb{Z} \} \cup
 \{ \f{\ell}{n (1 + \vep)} \in (a, b) : \ell \in \bb{Z} \} $.
Then, both $g(\lm)$ and $h (\lm)$ are constants for any $\lm \in (\al, \ba)$.
 \eeL

 \bpf
Since $\al$ and $\ba$ are two consecutive elements of the
ascending arrangement of all distinct elements of the set, it must
be true that there is no integer $\ell$ such that {\small $\al <
\f{\ell}{n (1 - \vep)}  < \ba$} or {\small $\al < \f{\ell}{n (1 +
\vep)} < \ba$}. It follows that there exist two integers $\ell$
and $\ell^\prime$ such that {\small $(\al, \ba) \subseteq \li (
\f{\ell}{n (1 - \vep)},  \f{\ell + 1}{n (1 - \vep)} \ri )$} and
{\small $(\al, \ba) \subseteq \li ( \f{\ell^\prime}{n (1 + \vep)},
\f{\ell^\prime + 1}{n (1 + \vep)} \ri )$}. Applying Lemma
\ref{minus_rev} and Lemma \ref{plus_rev}, we have {\small $g(\lm)
= g \li ( \f{\ell}{n (1 - \vep)} \ri )$} and {\small $h(\lm) = h
\li (\f{\ell^\prime + 1}{n (1 + \vep)} \ri )$} for any $\lm \in
(\al, \ba)$.

 \epf

\beL \la{lem_lim_rev}
 For any $\lm \in (0,1)$, $\lim_{\eta \downarrow 0} C(\lm + \eta) \geq C(\lm)$
 and $\lim_{\eta \downarrow 0} C(\lm - \eta) \geq C(\lm)$.
\eeL

\bpf

Observing that $h(\lm + \eta) \geq h(\lm)$ for any $\eta > 0$ and
that
 \bee g(\lm + \eta) & =  &  \lf n (\lm + \eta) (1 - \vep)
\rf + 1  \\
& = &  \lf n \lm (1 - \vep ) \rf + 1 + \lf n \lm (1 - \vep ) -
\lf n \lm ( 1  - \vep ) \rf + n \eta (1 - \vep) \rf \\
 & = & \lf n \lm ( 1 - \vep )
\rf + 1   = g(\lm) \eee for $0 < \eta < \f{ 1 + \lf n \lm( 1 - \vep
) \rf -  n \lm( 1 - \vep )} {n (1 - \vep)}$, we have \be
\la{ineqa_rev} S(n, g(\lm + \eta), h (\lm + \eta), \lm + \eta ) \geq
S(n, g(\lm), h (\lm), \lm + \eta ) \ee for $0 < \eta < \f{ 1 + \lf n
\lm( 1 - \vep ) \rf - n \lm( 1 - \vep )} {n (1 - \vep)}$.  Since
 \bee h(\lm + \eta) & =  &  \lc n (\lm + \eta) (1 + \vep)
\rc - 1\\
& = &  \lc n \lm (1 + \vep) \rc - 1 + \lc n \lm (1 + \vep) - \lc n
\lm (1 + \vep) \rc + n \eta (1 + \vep) \rc, \eee we have {\small \[
h(\lm + \eta) = \bec \lc n \lm(1 + \vep) \rc & \tx{for} \; n \lm( 1
+ \vep) = \lc n \lm( 1  + \vep) \rc \; \tx{and} \; 0 < \eta <
\f{1}{n (1 + \vep)},\\
\lc n \lm( 1 + \vep) \rc  - 1 & \tx{for} \; n \lm( 1 + \vep) \neq
\lc n \lm( 1 + \vep) \rc \; \tx{and} \; 0 < \eta < \f{\lc n \lm( 1 +
\vep) \rc - n \lm( 1 + \vep)}{n (1 + \vep)}. \eec
\]}
It follows that both $g(\lm + \eta)$ and $h(\lm + \eta)$ are
independent of $\eta$ if $\eta > 0$ is small enough. Since $S(n, g,
h, \lm + \eta)$ is continuous with respect to $\eta$ for fixed $g$
and $h$, we have that $\lim_{\eta \downarrow 0} S(n, g(\lm + \eta),
h (\lm + \eta), \lm + \eta )$ exists.  As a result, \bee \lim_{\eta
\downarrow 0} C(\lm + \eta) & = & \lim_{\eta \downarrow 0} S(n,
g(\lm + \eta), h (\lm + \eta), \lm + \eta )\\
& \geq & \lim_{\eta \downarrow 0} S(n, g(\lm), h (\lm), \lm + \eta )
= S(n,  g(\lm), h (\lm), \lm ) = C(\lm), \eee where the inequality
follows from (\ref{ineqa_rev}).

Observing that $g(\lm - \eta) \leq g(\lm)$ for any $\eta > 0$ and
that
 \bee h(\lm - \eta) & =  & \lc n( \lm - \eta) (1 + \vep)
\rc - 1\\
& = & \lc n \lm( 1 + \vep) \rc - 1 + \lc n \lm( 1 + \vep) -
\lc n \lm( 1  + \vep) \rc - n \eta ( 1 + \vep) \rc\\
 & = & \lc n \lm( 1 + \vep)
\rc - 1  = h(\lm) \eee for $0 < \eta < \f{ 1 + n \lm( 1 + \vep) -
\lc n \lm( 1 + \vep) \rc } {n ( 1 + \vep)}$, we have \be
\la{ineqb_rev} S(n, g(\lm - \eta), h (\lm - \eta), \lm - \eta ) \geq
S(n, g(\lm), h (\lm), \lm - \eta ) \ee for {\small $0 < \eta < \min
\li \{ \lm, \f{ 1 + n \lm( 1 + \vep) - \lc n \lm( 1 + \vep) \rc } {n
( 1 + \vep)} \ri \}$}. Since \bee g(\lm - \eta) & =  & \lf n( \lm -
\eta)(1 - \vep)
\rf + 1 \\
& = & \lf n \lm( 1 - \vep) \rf + 1 + \lf n \lm( 1 - \vep) - \lf n
\lm( 1 - \vep) \rf - n \eta (1 - \vep) \rf, \eee  we have {\small \[
g(\lm - \eta) = \bec \lf n \lm( 1 - \vep) \rf & \tx{for} \; n \lm( 1
- \vep) = \lf n \lm( 1  - \vep) \rf \; \tx{and} \; 0 < \eta <
\f{1}{n (1 - \vep)},\\
\lf n \lm( 1 - \vep) \rf  + 1 & \tx{for} \; n \lm( 1 - \vep) \neq
\lf n \lm( 1 - \vep) \rf \; \tx{and} \; 0 < \eta < \f{n \lm( 1 -
\vep) - \lf n \lm( 1 - \vep) \rf }{n (1 - \vep)}. \eec
\]}
It follows that both $g(\lm - \eta)$ and $h(\lm - \eta)$ are
independent of $\eta$ if $\eta > 0$ is small enough. Since $S(n, g,
h, \lm - \eta)$ is continuous with respect to $\eta$ for fixed $g$
and $h$, we have that $\lim_{\eta \downarrow 0} S(n, g(\lm - \eta),
h (\lm - \eta), \lm - \eta )$ exists. Hence, \bee \lim_{\eta
\downarrow 0} C(\lm - \eta) & =
& \lim_{\eta \downarrow 0} S(n, g(\lm - \eta), h (\lm - \eta), \lm -\eta )\\
& \geq & \lim_{\eta \downarrow 0} S(n, g(\lm), h (\lm), \lm -\eta )
= S(n, g(\lm), h (\lm), \lm ) = C(\lm), \eee where the inequality
follows from (\ref{ineqb_rev}).

\epf

By a similar argument as that of Lemma \ref{inbetween}, we have \beL
\la{inbetween_rev} Let $\al < \ba$ be two consecutive elements of
the ascending arrangement of all distinct elements of $\{a, b \}
\cup
 \{ \f{\ell}{n (1 - \vep)}  \in (a, b) : \ell \in \bb{Z} \} \cup
 \{ \f{\ell}{n (1 + \vep)} \in (a, b) : \ell \in \bb{Z} \} $.  Then,
 $C(\lm) \geq \min \{ C(\al), \; C(\ba) \}$ for any $\lm \in (\al, \ba)$.
\eeL

\bsk

Finally, to show Theorem \ref{thm_rev}, note that the statement about the coverage probability  follows
immediately from Lemma \ref{inbetween_rev}.  The number of elements of the finite set can be calculated by using
the property of the ceiling and floor functions.

\end{document}